\documentclass[oneside,10pt]{article}          
\usepackage[b5paper]{geometry}	    
\usepackage{graphicx,bm,epsf,amsfonts,amsmath,latexsym,amssymb} 
\usepackage{theorem}                
\usepackage{mathrsfs,upref}         
\usepackage{mathptmx}		    
\usepackage{jca}	            

\newtheorem{theorem}{Theorem}

\theoremstyle{definition}

\newtheorem{remark}{Remark}

\begin{document}
\def\f#1#2{\mbox{${\textstyle \frac{#1}{#2}}$}}
\def\dfrac#1#2{\displaystyle{\frac{#1}{#2}}}
\newcommand{\fr}{\frac{1}{2}}
\newcommand{\fs}{\f{1}{2}}
\newcommand{\g}{\Gamma}
\newcommand{\br}{\biggr}
\newcommand{\bl}{\biggl}
\newcommand{\ra}{\rightarrow}

\title[Zeta function integrals]{Integrals of products of Hurwitz zeta functions and
                    the Casimir effect in $\bm{\phi^4}$ field theories}


\author{M. A. Shpot, M. P. Chaudhary and R. B. Paris}

\address{M. A. Shpot, Institute for Condensed Matter Physics, 79011 Lviv, Ukraine\\
\email{shpot.mykola@gmail.com}}

\address{M. P. Chaudhary, Department of Mathematics, International Scientific Research
and Welfare Organization, New Delhi, 110018 India\\
\email{dr.m.p.chaudhary@gmail.com}}

\address{R. B. Paris, Division of Computing and Mathematics, University of Abertay Dundee,
Dundee DD1 1HG, UK\\
\email{r.paris@abertay.ac.uk}}

\CorrespondingAuthor{R. B. Paris}

\date{12.06.2016}                               



\thanks{The authors wish to acknowledge M. Katsurada
for pointing out the result stated in (\ref{e35}).
MAS and MPC are grateful to Professors H.M. Srivastava and J. Choi for
giving their opinion on an early draft of the paper.
Part of the work of MAS was carried out at
the Fakult\"at f\"ur Physik of the Universit\"at Duisburg-Essen and he is grateful to
Prof. H.W. Diehl for his warm hospitality and financial support
}

\begin{abstract}
We evaluate two integrals over $x\in [0,1]$ involving products of the function
$\zeta_1(a,x)\equiv \zeta(a,x)-x^{-a}$ for $\Re (a)>1$, where $\zeta(a,x)$ is the Hurwitz zeta function.
The evaluation of these integrals for
the particular case of integer $a\geq 2$ is also presented.
As an application we calculate the $O(g)$ weak-coupling expansion
coefficient $c_{1}(\varepsilon)$ of the Casimir
energy for a film with Dirichlet-Dirichlet boundary conditions, first stated by Symanzik
[\textit{Schr{\"o}dinger representation and Casimir effect in renormalizable quantum field theory,}
Nucl. Phys. B {\bf 190} (1981) 1-44] in the framework of $g\phi^4_{4-\varepsilon}$ theory.
\end{abstract}

\maketitle



\section{Introduction}

The Hurwitz zeta function $\zeta(a,x)$, 
see (\ref{ZD}), is one of the most fundamental functions in mathematics.
It has important applications, for example, in number theory \cite{Apostol13,Arakawa14},
probability theory \cite{Karatsuba92} and in the evaluation of mathematical constants
\cite{Finch03}.
This function also finds application in numerous areas of mathematical physics.
Examples are in the evaluation of functional determinants and
in gases of free electrons in the presence of a magnetic field,
calculations of the Casimir effect \cite{Cas48} in quantum field theory and
its analogues in statistical physics and thermal field theory at high temperatures;
see e.g. \cite{Krech,Elizalde12,KG99,KR04} for extensive lists of physical examples.

The Hurwitz zeta function is defined by the series \cite[p. 89]{Hurwitz1882}, \cite[Ch. 13]{WhW}
\begin{equation}\label{ZD}
\zeta(a,x)=\sum_{k=0}^\infty\frac{1}{(k+x)^a}\qquad(\Re(a)>1;\ x\neq 0, -1, -2, \ldots)
\end{equation}
and elsewhere by analytic continuation,
apart from $a=1$, where it has a simple pole with unit residue (see (\ref{LE}) below).
In its convergence domain, the series (\ref{ZD}) converges absolutely and uniformly.
The function $\zeta(a,x)$ reduces to the Riemann zeta function when $x=1$, viz.
\[\zeta(a,1)=\zeta(a)=\sum_{k=1}^\infty\frac{1}{k^a}\qquad(\Re(a)>1),\]
and as $x\to0^+$, its behavior is singular described by $\zeta(a,x)\sim x^{-a}$.

Before proceeding, we record some necessary preliminary results related to $\zeta(a,x)$.
By a straightforward manipulation of the above series definitions it is easy to obtain
the identity
\begin{equation}\label{e13}
\zeta(a,2)=\zeta(a)-1\,,
\end{equation}
which will be used below.
It is worth quoting here the Wilton formula (see, for example, \cite[Eq.~(8)]{KKY01} or
\cite[p.~248, Eq.~(7)]{SC} and references therein)
\begin{equation}\label{TE}
\zeta(a,b-z)=\sum_{k=0}^\infty\frac{(a)_k}{k!}\,\zeta(a+k,b)\,z^k\,\qquad (|z|<|b|,\;a\ne1)\,,
\end{equation}
where $(a)_k$ denotes the Pochhammer symbol
\begin{equation}\label{PH}
(a)_k=\frac{\g(a+k)}{\g(a)}=a(a+1)\ldots(a+k-1),\qquad (0)_0 = 1.
\end{equation}
The relation (\ref{TE}) can be regarded as a Taylor
expansion of the function on the left-hand side in powers of
$z$ and will be used below in the special cases $b=1$ and $b=2$.
At $a=1$, which is forbidden in (\ref{TE}), the poles of the Hurwitz functions cancel
on both sides due to the Laurent expansion \cite[p. 159, Eq. (15)]{SC}
\begin{equation}\label{LE}
\zeta(1+\epsilon,z)=\frac1\epsilon\Big\{1-\epsilon\psi(z)+O(\epsilon^2)\Big\}\,.
\end{equation}
This transforms the formula (\ref{TE}) to
\begin{equation}
\psi(b-z)=\psi(b)-\sum_{k=1}^\infty\zeta(1+k,b)\,z^k\,\qquad (|z|<|b|),
\end{equation}
where $\psi(b)=\g'(b)/\g(b)$ is the psi-function or logarithmic derivative of the
Euler gamma function.

Integrals of Hurwitz zeta functions with respect to the variable $x$  have been discussed by several authors
\cite{KL52,KKY01,EM02,EM022}.
In \cite{EM02}\footnote{
This reference contains a short review of availability of such
integrals in standard tables of integrals.}
(see also \cite{LHK09}), the Hurwitz transform
\[\int_0^1 f(x) \zeta(a,x)\,dx\]
was considered for a variety of functions $f(x)$, including integer powers of $x$,
trigonometric, exponential, and other functions.
In \cite{KL52}, the evaluation of the mean value integral $\int_0^1|\zeta_1(a,x)|^2 dx$
was initiated. This problem has been considered further in
\cite{And92,KM93,KM96,KTY06} and papers cited therein,
in a number-theoretic context when $a=\fs+it$ and $t\ra\pm\infty$.
In this last integral, the auxiliary zeta function $\zeta_1(a,x)$ appears, obtained by
splitting the singular zeroth term from $\zeta(a,x)$:
\begin{equation}\label{e11}
\zeta_1(a,x)=\zeta(a,x)-x^{-a}.
\end{equation}

The advantage of the function \cite{KL52}
\begin{equation}\label{e12}
\zeta_1(a,x)=\sum_{k=1}^\infty\frac{1}{(k+x)^a}=\zeta(a,x+1)\qquad(\Re(a)>1)
\end{equation}
is (cf. (\ref{ZD})) that it is continuous in the whole interval $x\in\,[0,1]$. Note that $\zeta_1(a,0)=\zeta(a)$.
It is easy to deduce from the Wilton formula (\ref{TE}) an analogous Taylor expansion
for $\zeta_1(a,b-z)$ by using its definition in (\ref{e11}) and employing the
binomial expansion in the term $(b-z)^{-a}$.
In the special case of $b=1$, to be considered below, this reads
\begin{equation}\label{WD}
\zeta_1(a,1-z)=\sum_{k=0}^\infty\frac{(a)_k}{k!}\,
\{\zeta(a+k)-1\}\,z^k\,\qquad (|z|\leq 1,\;a\ne1)\,.
\end{equation}
Note that the last formula has a well-defined limit as $z\to 1$ with the left-hand side reducing to $\zeta(a)$, whereas both sides of (\ref{TE}) diverge at $z=b$.

We shall also use the integral representation for the function $\zeta_1(a,x)$ given by
\begin{equation}\label{e23}
\zeta_1(a,x)=\frac{1}{\Gamma(a)}\int_0^\infty \frac{t^{a-1} e^{-xt}}{e^t-1}\,dt\,\qquad(\Re(a)>1)\,,
\end{equation}
which follows from the standard formula \cite[p. 90]{Hurwitz1882}
(or, for example, \cite[p. 155, Eq. (2)]{SC} with $a\mapsto x+1$).

In the present article we obtain expressions for the integrals
involving products of two auxiliary zeta functions $\zeta_1$, namely
\[\int_0^1 \zeta_1(a,x) \zeta_1(b,x)\,dx \qquad\mbox{and}\qquad \int_0^1 \zeta_1(a,x) \zeta_1(b,1-x)\,dx\]
for $\Re (a, b)>1$. The first integral has been given in \cite[Eq.~(5)]{And92}%
\footnote{The sign at the second term on the right-hand side of \cite[Eq.~(5)]{And92} is misprinted.
The correct contribution is $-1/(1-z-w)$.}
and \cite[Corollary 4]{KM96};
here we give a different derivation. Both evaluations of the above integrals present
apparent singularities at integer values of $a, b\geq 2$ and require a limiting
procedure to determine their values for these special values of the parameters.

As an application of the above two integrals, we give a derivation of the Casimir amplitude
resulting from the attraction of two infinite parallel planes with Dirichlet boundary conditions.
This result was stated without proof in \cite[p.~12]{Sym81}
and has been derived subsequently by other means in \cite{KD92}.

Finally, we note that integrals of the same form have been calculated very recently in
\cite{HKK15}.
They involve the alternating counterpart \cite[Eq. (1.5)]{CvSr09} of the Hurwitz zeta function,
\[
\zeta_E(a,x):=\sum_{k=0}^\infty\frac{(-1)^k}{(k+x)^a}\qquad(\Re(a)>0;\ x\neq 0, -1, -2, \ldots).
\]
Its relation to the usual $\zeta(a,x)$ is \cite[Eq. (1.6)]{HKK15}, \cite[Eq. (1.6)]{CvSr09}
\[
\zeta_E(a,x)=2^{1-a}\,\zeta\Big(a,\frac x2\Big)-\zeta(a,x)
=2^{-a}\Big[\zeta\Big(a,\frac x2\Big)-\zeta\Big(a,\frac{x+1}2\Big)\Big].
\]

\section{The integral of a product of auxiliary Hurwitz zeta functions}

We first consider the integral
\begin{equation}\label{e21}
I(a,b):=\int_0^1 \zeta_1(a,x) \zeta_1(b,x)\,dx,
\end{equation}
where the parameters $a$ and $b$ satisfy $\Re (a, b)>1$.
This region comes from the condition that the simple series definition (\ref{e12})
of the functions $\zeta_1(a,x)$ and $\zeta_1(b,x)$ exists, and there is no additional constraint
due to the convergence of the integral in (\ref{e21}). The integrand is a smooth
integrable function everywhere in the integration region $[0,1]$.
This integral has been evaluated in \cite[Eq.~(5)]{And92}
(see footnote 2) and \cite[Corollary 4]{KM96}
with $a, b\neq 2, 3 \ldots\ $.
The main motivation for taking different parameters $a$ and $b$
in the integrand of (\ref{e21}) has been in evaluation of the "mean-square" of the Hurwitz zeta function, when $a$ and $b$ are complex conjugate numbers,
especially, $a=b^*=\fs+it$ (see \cite{KL52,And92,KM93,KM96,KTY06}).
In view of the importance of this problem and for
completeness in presentation, we give here a different method of proof.

As in \cite{KM96} (cf. \cite[Eq. (3.1)]{KM96}), we start by using Atkinson's dissection \cite{Atk49} to
split the double summation in
$\zeta_1(a,x) \zeta_1(b,x)=\sum_{k,n\ge 1}(\cdot)$
into parts with $k=n$, $k<n$, and $k>n$, which yields
\[\zeta_1(a,x)\zeta_1(b,x)=\zeta_1(a+b,x)+\sum_{n=1}^\infty \left\{\frac{\zeta_1(a,n+x)}{(n+x)^{b}}+\frac{\zeta_1(b,n+x)}{(n+x)^{a}}\right\}\qquad(\Re (a, b)>1).\]
Thus we find
\begin{equation}\label{e22}
\int_0^1\zeta_1(a,x)\zeta_1(b,x)\,dx=\int_0^1\zeta_1(a+b,x)\,dx+\sum_{n=1}^\infty\int_0^1
\left\{\frac{\zeta_1(a,n+x)}{(n+x)^{b}}+\frac{\zeta_1(b,n+x)}{(n+x)^{a}}\right\}dx,
\end{equation}
where the first term on the right-hand side is straightforward and given by
\begin{equation}\label{e22a}
\int_0^1\zeta_1(a+b,x)\,dx=\frac{1}{a+b-1}\qquad (\Re(a+b)>1).
\end{equation}

The second and third terms of (\ref{e22}) are similar and so it suffices to present the details of only one of these terms. Choosing to evaluate the third term, we can use the integral representation (\ref{e23})
for the function $\zeta_1$ and the Laplace integral
\[(n+x)^{-a}=\frac{1}{\Gamma(a)}\int_0^\infty y^{a-1} e^{-(n+x) y}\,dy \qquad(\Re(a)>1).\]
Substitution of these results then yields
\begin{align*}
S:=&\sum_{n=1}^\infty\int_0^1 \frac{\zeta_1(b,n+x)}{(n+x)^a}\,dx
\\
=&\frac{1}{\Gamma(a)\Gamma(b)}\sum_{n=1}^\infty\int_0^\infty \, y^{a-1}\int_0^\infty
t^{b-1}\frac{e^{-n(t+y)}}{e^t-1}\bl\{\int_0^1 e^{-x(t+y)}\,dx\br\}\,dt\,dy.
\end{align*}

Making use of the evaluations
\[\sum_{n=1}^\infty e^{-n(t+y)}=\frac{1}{e^{t+y}-1}=\frac{e^{-(t+y)}}{1-e^{-(t+y)}}
\qquad\mbox{and}\qquad \int_0^1 e^{-x(t+y)}\,dx
=\frac{1-e^{-(t+y)}}{t+y},\]
and noting that the factors $1-e^{-(t+y)}$ cancel\footnote{It is actually this cancellation that allowed us to easily finish the calculation.
There has been no such cancellation in our attempts to calculate the integral
$J(a,b)$,
considered in the next section, by employing Atkinson's dissection in the product
$\zeta_1(a,x)\zeta_1(b,1-x)$ in its integrand. For this reason we had to resort
to alternative means in the evaluation of $J(a,b)$
in Section 3.} in the product of these two results,
we can write the above expression as
\begin{equation}\label{A25}
S=\frac{1}{\Gamma(a)\Gamma(b)}\int_0^\infty \frac{t^{b-1}e^{-t}}{e^t-1}
\bl\{\int_0^\infty \frac{y^{a-1}e^{-y}}{t+y}\,dy\br\}\,dt
=\frac{1}{\Gamma(b)}\int_0^\infty\frac{t^{a+b-2}}{e^t-1}\,\Gamma(1-a,t)\,dt,
\end{equation}
where the integral over $y$ has been evaluated in terms of the incomplete gamma function given by \cite[Eq.~(8.6.4)]{DLMF}
\[\g(1-a,t)=\frac{t^{1-a}e^{-t}}{\g(a)} \int_0^\infty \frac{y^{a-1} e^{-y}}{t+y}\,dy \qquad (\Re (a)>0).\]
Now we take into account that
\[\Gamma(1-a,t)=\Gamma(1-a)-\gamma(1-a,t)\]
where $\gamma(a,x)$ denotes the lower incomplete gamma function.
Performing the integral over $t$ in (\ref{A25}) related to the first term
of this last equation, we obtain by (\ref{e23})
\[S=B(a+b-1,1-a)\,\zeta(a+b-1)-
\frac{1}{\Gamma(b)}\int_0^\infty \frac{t^{a+b-2}}{e^t-1}\,\gamma(1-a,t)\,dt,\]
where $B(x,y)$ is the beta function
\begin{equation}\label{BF}
B(x,y)=\frac{\g(x) \g(y)}{\g(x+y)}\,.
\end{equation}

In the remaining integral over $t$ we use the series expansion \cite[Eq.~(8.5.1)]{DLMF}
\[\gamma(1-a,t)=t^{1-a}e^{-t}\sum_{n=0}^\infty\frac{t^n}{(1-a)_{n+1}}\qquad(|t|<\infty),\]
where again $(a)_k$ is the Pochhammer symbol (see (\ref{PH})). So,
\begin{align*}
S&=B(a+b-1,1-a)\,\zeta(a+b-1)-
\frac{1}{\Gamma(b)}\sum_{n=0}^\infty\frac{1}{(1-a)_{n+1}}
\int_0^\infty \frac{t^{b+n-1}e^{-t}}{e^t-1}\,dt\\
&=B(a+b-1,1-a)\,\zeta(a+b-1)-
\sum_{n=0} ^\infty \frac{(b)_n}{(1-a)_{n+1}}\,\zeta(b+n,2)
\end{align*}
by application of (\ref{e23}) again.
Finally this leads us, via (\ref{e13}) and (\ref{e22}), to the following theorem:

\begin{theorem} For $\Re (a, b)>1$, we have the evaluation
\[\int_0^1 \!\zeta_1(a,x) \zeta_1(b,x)dx
=\frac{1}{a{+}b{-}1}+\bl\{\!B(a+b-1,1-a)+B(a+b-1,1-b)\!\br\}\zeta(a+b-1)\]
\begin{equation}\label{e25}
\hspace{2cm}-\sum_{n=0}^\infty\frac{(a)_n}{(1-b)_{n+1}}\{\zeta(a+n)-1\}
-\sum_{n=0}^\infty\frac{(b)_n}{(1-a)_{n+1}}\{\zeta(b+n)-1\}.
\end{equation}
\end{theorem}
Though derived for $\Re(a, b)>1$, this result provides the analytic continuation
of $I(a,b)$ into the whole complex $a$, $b$ planes except at $a=1$ and $b=1$
since both sides of (\ref{e25}) are analytic functions of $a$ and $b$,
apart from $a=1$, $b=1$ where there is a double pole.

\begin{remark} The second term in (\ref{e25})
coincides with the
similar Mikol\'as integral of the product of \emph{Hurwitz} zeta functions
(cf. \cite[Eq. (6.4)]{KM96}, \cite[Eq. (2)]{And92}, \cite[Eq. (3.2)]{EM02})
\[
\int_0^1\zeta(a,x)\zeta(b,x)\,dx=
\left\{B(a+b-1,1-a)+B(a+b-1,1-b)\right\}\zeta(a+b-1)
\]
valid for $\Re(a), \Re(b), \Re(a+b)<1$.
\end{remark}

The expression on the right-hand side of (\ref{e25}) presents apparent singularities when
$a, b=2, 3, \ldots$ due to the presence of the factors $\g(1-a)$, $(1-a)_{n+1}$ and $\g(1-b)$, $(1-b)_{n+1}$ in the last three terms.
However, when $a=b$
these singularities cancel to leave a regular function.
The value of $I(a,b)$ when $a=b=m$ at integer values of $m\geq 2$
is obtained by a limiting process and is considered in the appendix; see also \cite{KM97}.
This yields
\begin{theorem} For integer $m=2, 3, \ldots\ $, we have the evaluation
\begin{align}\label{e28}
\int_0^1 \zeta_1^2(m,x)\,dx=\frac{1}{2m-1}&+\frac{2(-)^m}{\g^2(m)}\bl\{\g(2m-1) Z(m)
\\\nonumber
&-\sum_{n=m-1}^\infty \frac{(n\!+\!m\!-\!1)!}{(n\!-\!m\!+\!1)!}\,\psi(n+2-m)\,\{\zeta(n+m)-1\}\br\}
\\\nonumber
&+\frac{2}{\g^2(m)}\sum_{n=0}^{m-2}(-)^n\,(n\!+\!m\!-\!1)!\,(m\!-\!n\!-\!2)!\,\{\zeta(n+m)-1\},
\end{align}
where
\[Z(m):=\psi(2m-1) \zeta(2m-1)+\zeta'(2m-1)\]
and the prime on the zeta function denotes its derivative.
\end{theorem}
It is readily shown that this result is equivalent to that obtained in \cite[Theorem 3]{KM97}.

\section{The integral of a product of auxiliary Hurwitz zeta functions with complementary arguments}

The second integral we consider is
\begin{equation}\label{e31}
J(a,b):=\int_0^1\zeta_1(a,x)\zeta_1(b,1-x)\,dx \qquad (a, b\neq 1).
\end{equation}
As in the case of $I(a,b)$ in (\ref{e21}),
the integrand in $J(a,b)$ is regular in the whole integration range $x\in[0,1]$ and the integral
converges for any allowed values of $a$ and $b$. Following the approach of Andersson \cite{And92},
we first evaluate the integral (\ref{e31}) for $0<\Re (a, b)<1$
and then extend it into $\Re (a, b)>1$ by analytic continuation.
The result we obtain again presents apparent singularities at
$a,\, b=2, 3, \ldots\ $; the value of $J(a,b)$ in the special case $a=b$ for integer values of
$a$ is derived by a limiting process similar to that in the appendix.

\subsection{Evaluation of $\bm{J(a,b)}$}

In a similar manner to that employed in \cite{And92}, we rewrite the integral $J(a,b)$ as follows:
\begin{align}
J(a,b)&=\int_0^1\left[\zeta(a,x)-x^{-a}\right]\left[\zeta(b,1-x)-(1-x)^{-b}\right] dx\nonumber\\
&=\int_0^1\zeta(a,x)\zeta(b,1-x)\,dx-\int_0^1 \{x^{-a}\zeta_1(b,1-x)+x^{-b} \zeta_1(a,1-x)\}\,dx\nonumber\\
&\hspace{6cm}-\int_0^1 x^{-a}(1-x)^{-b}\,dx.\label{e32}
\end{align}
An explicit expression for the first integral on the right-hand side of (\ref{e32})
can be obtained from the result derived in \cite[Eq.~(3.4)]{EM02} given by
\begin{equation}\label{e33}
\int_0^1\zeta(a,x)\zeta(b,1-x)\,dx=B(1-a,1-b)\zeta(a+b-1)
\end{equation}
for the region $\Re (a)\leq 0$, $\Re (b)\leq 0$,
where $B(x,y)$ is the beta function defined in (\ref{BF}).
The right-hand side of (\ref{e33}) can be extended to the domain $\Re (a)<1$, $\Re (b)<1$
by analytic continuation; hence (\ref{e33}) is valid for $\Re (a)<1$ and $\Re (b)<1$.

To calculate the second integral in (\ref{e32}) we use
the Taylor expansion (\ref{WD})
and perform a term-by-term integration. This gives (with $0<\Re (a, b)<1$)
\begin{equation}\label{e34}
\int_0^1 x^{-a}\zeta_1(b,1-x)\,dx=\sum_{n=0}^\infty\frac{(b)_n}{n!}\,\frac{\zeta(b+n)-1}{1-a+n}\,.
\end{equation}
A similar result holds for the third integral in (\ref{e32}) with $a$ and $b$ interchanged.

Finally, evaluation of the last integral on the right-hand side of
(\ref{e32}) as $B(1-a,1-b)$ (when $0<\Re (a, b)<1$) then yields the result:

\begin{theorem} When $0<\Re (a, b)<1$, we have the evaluation
\[\int_0^1\zeta_1(a,x)\zeta_1(b,1-x)\,dx =B(1-a,1-b)\{\zeta(a+b-1)-1\}\hspace{3cm}\]
\begin{equation}\label{e35}
\hspace{3.2cm}-\sum_{n=0}^\infty\frac{(a)_n}{n!}\,\frac{\zeta(a+n)-1}{n+1-b}
-\sum_{n=0}^\infty\frac{(b)_n}{n!}\,\frac{\zeta(b+n)-1}{n+1-a}~.
\end{equation}
Although calculated for $0<\Re (a, b)<1$, this result provides the analytic continuation of $J(a,b)$ into $\Re (a, b)>1$
since both sides of (\ref{e35}) are analytic functions of $a$ and $b$,
except at $a=1$, $b=1$ where there is a double pole.
\end{theorem}

Again, the expression on the right-hand side of (\ref{e35}) presents
apparent singularities at $a, b=2, 3, \ldots$, but these will be
shown to cancel when $b=a$ to leave a regular function of $a$ when $\Re (a)>1$.

\begin{remark} The integral in (\ref{e34}) is the counterpart of the similar integral
\begin{equation}\label{e34b}
\int_0^1 x^{-a} \zeta_1(b,x)\,dx=\sum_{n=0}^\infty \frac{(b)_n}{(1-a)_{n+1}}\,\{\zeta(b+n)-1\}
\qquad (\Re (a)<1,\ b\ne1, 0, -1, -2,...)
\end{equation}
obtained in \cite[Eq.~(4)]{And92}.
\end{remark}

\subsection{Evaluation of $\bm{J(a,a)}$ for integer $\bm{a\geq 2}$}

We consider the special case of $J(a,b)$ defined in (\ref{e31}) when $a=b$.
We let $a=m+\epsilon$,
where $m=2, 3, \ldots\,$, and consider the limit $\epsilon\ra 0$.
Using the definition of the beta
function (\ref{BF}), we find from (\ref{e35})
\begin{align}\label{JMM}
J&(m,m)=\\\nonumber
&\lim_{\epsilon\ra 0}\bl\{\frac{\g^2(1-m-\epsilon)}{\g(2-2m-2\epsilon)}\,
\{\zeta(2m-1+2\epsilon)-1\}-2\sum_{n=0}^\infty\frac{(m+\epsilon)_n}{n!}\;
\frac{\zeta(n+m+\epsilon)-1}{n+1-m-\epsilon}\br\}.
\end{align}
Extracting the singular term from the sum, we can then write
\begin{equation}\label{e37}
J(m,m)=\lim_{\epsilon\ra 0} D(m+\epsilon)-2\sum_{n=0}^\infty{}\!' \frac{(m)_n}{n!}\,\frac{\zeta(n+m)-1}{n+1-m},
\end{equation}
where the prime on the summation sign denotes the omission of the term corresponding to $n=m-1$ and we have defined
\[
D(m+\epsilon):=\frac{\g^2(1-m-\epsilon)}{\g(1-s-2\epsilon)}\,
\{\zeta(s+2\epsilon)-1\}+\frac{2\g(s+\epsilon)}
{\epsilon\,\g(m) \g(m+\epsilon)}\,\{\zeta(s+\epsilon)-1\}.
\]
In writing $D(m+\epsilon)$ in the above form we have used the definition (\ref{PH})
for the Pochhammer symbol $(m+\epsilon)_{m-1}$ appearing in the
$n=(m-1)$th term of the sum in (\ref{JMM}) and temporarily introduced
the shorthand $s:=2m-1\ge 3$.

When expanding $D(m+\epsilon)$ for $\epsilon\ll 1$,
we use the well-known small-parameter expansions of the Euler gamma function
\begin{equation}\label{GE}
\g(z+\epsilon)=\g(z)\{1+\epsilon \psi(z)+O(\epsilon^2)\}
\end{equation}
where $\psi(z)$ denotes again the psi-function, and, in particular,
\begin{equation}\label{EG}
\g(-p+\alpha)=\frac{(-)^p}{p!\,\alpha}\{1+\alpha\psi(p+1)+O(\alpha^2)\}
\end{equation}
if $p$ is a non-negative integer. We need these two formulas with
$z=s$, $z=m$, and $p=s-1$, $p=m-1$, respectively.
By a  straightforward calculation we see immediately that the pole contributions
\[
\pm\,\frac{2\g(s)}{\epsilon\,\g^2(m)}\,\{\zeta(s)-1\}
\]
of the two terms in $D(m+\epsilon)$ mutually cancel, and we are left with their
finite parts which combine to yield
\[
D(m+\epsilon)=
-\frac{2\g(s)}{\g^2(m)}\Big\{[\psi(s)-\psi(m)]\,[\zeta(s)-1]+\zeta'(s)+O(\epsilon)\Big\}.
\]
Hence we obtain, by recalling (\ref{e37}) and $s=2m-1$,
\begin{theorem} For $m=2, 3, \ldots\,$, we have the evaluation
\begin{align}\nonumber
\int_0^1\zeta_1(m,x)\zeta_1(m,1-x)\,dx=&
\\\nonumber
\frac{2\g(2m-1)}{\g^2(m)}
&\Big\{\left[\psi(m)-\psi(2m-1)\right]\left[\zeta(2m-1)-1\right]-\zeta'(2m-1)\Big\}
\\\label{e38}
&-2\sum_{n=0}^\infty{}\!' \frac{(m)_n}{n!}\,\frac{\zeta(n+m)-1}{n+1-m},
\end{align}
where the prime on the summation denotes the omission of the term corresponding to $n=m-1$.
\end{theorem}

\section{An application: The Casimir amplitude up to O($\bm g$) in the $\bm{g\phi^4}$ field theory}

In 1981, Symanzik \cite{Sym81} proved the finiteness of the
Casimir energy%
\footnote{This is a counterpart of the interaction energy between two
parallel conducting plates due to zero-point energy fluctuations
of the electromagnetic field calculated by Casimir in 1948 \cite{Cas48}.
For more extended information and various applications the interested reader
may be referred to, for example, \cite{Krech,KG99,bordag09}.}
$E(g,L,\mu,\varepsilon)$ of a pair of parallel plates with Dirichlet boundary
conditions to all orders in the perturbation expansion of the massless
$g\mu^\varepsilon\phi^4_d$ theory.
He considered the (super) renormalizable Euclidean scalar $\phi^4$ field theory
with the dimensionless coupling constant $g$ in $d=4-\varepsilon$
space dimensions for $\varepsilon\in[0,1)$.
The $d-1$ dimensional parallel boundary plates are separated by the distance $L$.
An arbitrary momentum scale $\mu$ of dimension [length]$^{-1}$ is introduced into the theory;
it combines, for example, with the slab thickness $L$ to form the dimensionless
combination $\mu L$.
The position vector in the slab $\bm x=(\bm r,z)$ consists of two components:
the vector $\bm r\in\mathbb R^{d-1}$ is parallel to the boundary planes, whereas
$z\in[0,L]$ is the one-dimensional coordinate running in the direction perpendicular to them.

The Casimir energy density has been given in \cite[Eq. (4.5b)]{Sym81} as
an infinite series
\begin{equation}\label{CE}
E(g,L,\mu,\varepsilon)/A=L^{-(d-1)}\Big[
c_0(\varepsilon)+c_1(\varepsilon)\bar g(g,\mu L,\varepsilon)+O(\bar g^2)\Big]
\end{equation}
in powers of the running coupling constant $\bar g=\bar g(g,\mu L,\varepsilon)$
with computable coefficients $c_n(\varepsilon)$, where $A$ is the area of the boundary planes.
The coefficient $c_0(\varepsilon)$ corresponds to the one-loop contribution
of the perturbation theory given by the Feynman diagram
$\;\;\raisebox{-0.9pt}{\includegraphics[width=10pt]{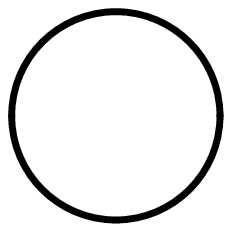}}$.
The value of $c_0(\varepsilon)$ can be found in \cite[Eq. (4.4)]{Sym81},
\[
c_0(\varepsilon)=-(4\pi)^{-d/2}\Gamma(d/2)\zeta(d),
\]
along with a sketch of its derivation.
Its expansion up to the first order in $\varepsilon$ is
\begin{equation}\label{C0}
c_0(\varepsilon)=-\frac{\pi^2}{1440}\left\{1+\Big[\frac{\gamma-1}{2}
+\frac12\ln(4\pi)-\frac{\zeta'(4)}{\zeta(4)}\Big]
\varepsilon\right\}+O(\varepsilon^2)\,,
\end{equation}
where $\gamma=0.57721\ldots$ is Euler's constant.

The coefficient $c_1(\varepsilon)$ comes from the two-loop, $O(g)$ Feynman diagram 
$\raisebox{-1.0mm}{\includegraphics[width=8mm]{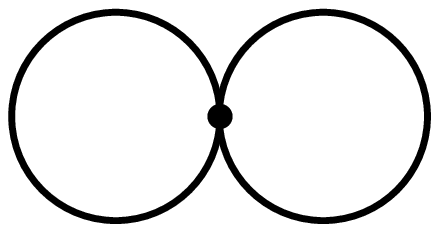}}\;$.
This is given by the expression
(see \cite[Eq.~(4.6)]{Sym81})
\[L^{-3+\varepsilon} c_1(\varepsilon) {\bar g}
(g,\mu L,\varepsilon)=
\f{1}{8}g\mu^\varepsilon\bl[\int_0^L\bl\{\left[G_D^L(0;z,z)-G_b(0)\right]^2
-\left[G_D(0;z,z)-G_b(0)\right]^2\]
\begin{equation}\label{e40}
-[G_D(0;L\!-\!z,L\!-\!z)-G_b(0)]^2\br\}dz
-2 \int_L^\infty \left[G_D(0;z,z)-G_b(0)\right]^2dz\br]+O(g^2)\,.
\end{equation}
To explain the functions appearing in the integrand it is useful to start with
the free-field propagator (Green function) of the massless $\phi^4_d$ theory
in the infinite $d$ dimensional space, $G_b(\bm x)$. It is given by
\[
G_b(x)=Cx^{-2+\varepsilon}
\qquad\mbox{with}\qquad C=\f{1}{4}\pi^{-2+\varepsilon/2}\Gamma(1-\fs\varepsilon),
\]
where $x$ is the distance between any two points in $\mathbb R^d$.
Its (infinite) value at $x=0$ is subtracted from each term of (\ref{e40}).
Further, the $G_D(\bm r;z_1,z_2)$ is the Dirichlet propagator in the semi-infinite space
confined by a single boundary plane. It is given by
\[
G_D(r;z_1,z_2)=G_b\left(r;|z_1-z_2|\right)-G_b\left(r;z_1+z_2\right),
\]
where $r=|\bm r_1-\bm r_2|$ is the "parallel" distance between the points
$\bm x_1$ and $\bm x_2$. Its explicit expression, though in
somewhat different notation, is given in \cite[Eq.~(2.13)]{Sym81}.
Finally, the Dirichlet-Dirichlet Green function in the film geometry
described at the beginning of the present section is given by
(compare with the previous formula)
\[
G^L_D(r;z_1,z_2)=\sum_{n=-\infty}^{\infty}\Big[
G_b\left(r;|z_1-z_2+2Ln|\right)-G_b\left(r;z_1+z_2+2Ln\right)\Big]\,.
\]
Again, the ``expanded" version of this formula can be found in \cite[Eq.~(4.7)]{Sym81}.

For coincident points $\bm x_1$ and $\bm x_2$ inside the slab, which is needed
to construct the $O(g)$ graph
$\raisebox{-1.0mm}{\includegraphics[width=8mm]{g2h}}\;\;$,
we have for the combination in the first square brackets in (\ref{e40})
\cite[Eq.~(4.8)]{Sym81}
\[G_D^L(0;z,z)-G_b(0)
=C\,(2L)^{-2+\varepsilon}\left[ 2 \zeta(2-\varepsilon)-\zeta(2-\varepsilon,z/L)
-\zeta(2-\varepsilon,1-z/L)\right],\]
where the functions $\zeta(a,x)$ are the Hurwitz zeta functions
defined in (\ref{ZD}).
Finally, for the second square brackets in (\ref{e40}) we have
\[G_D(0;z,z)-G_b(0)=-G_b(2z)=-C(2z)^{-2+\varepsilon}\]
along with a similar expression inside the third brackets.

After giving all these prerequisites for the calculation,
Symanzik \cite[p. 12]{Sym81} just gave the final result (quoted in (\ref{e48})
below) for $c_1(\varepsilon)$ without referring to any details of the calculation.
However, we find these details essential and consequently re-derive his result by
using just the above formulas given in purely coordinate representation.%
\footnote{Ten years after the paper by Symanzik \cite{Sym81}, Krech and Dietrich \cite{KD92} reproduced
the coefficient $c_1(\varepsilon)$ by using the momentum representation
along with dimensional and zeta regularizations.}

To proceed, we express the first term on the right-hand side of (\ref{e40}) as
\begin{equation}\label{e41}
\f{1}{8}g\mu^\varepsilon(I_{0,L}-2I_{L,\infty})C^2,
\end{equation}
where the related Feynman integrals $I_{0,L}$ and $I_{L,\infty}$
are labeled with respect to their integration ranges over $z$.
These are given by
\begin{equation}\label{e44}
I_{L,\infty}=\int_L^\infty (2z)^{-4+2 \varepsilon} dz=
2^{-4+2 \varepsilon}\frac{L^{-3+2 \varepsilon}}{3-2 \varepsilon}=
2^{-4+2 \varepsilon}L^{-3+2 \varepsilon}\frac{1}{2\alpha-1}\qquad (\alpha>\fs)
\end{equation}
and
\begin{equation}\label{e42}
I_{0,L}=2^{-4+2 \varepsilon} L^{-3+2 \varepsilon} K(\alpha),
\end{equation}
where
\begin{equation}\label{e43}
K(\alpha)=\int_0^1\left\{\left[2\zeta(\alpha)-\zeta(\alpha,x)-\zeta(\alpha,1-x)\right]^2
-x^{-2\alpha}-(1-x)^{-2\alpha}\right\}dx\quad(\alpha\in(1,2])
\end{equation}
with $x=z/L$.

For convenience we have introduced the notation $\alpha\equiv 2-\varepsilon$.
This will facilitate the comparison with expressions of the preceding sections.
Note that for $0\le\varepsilon< 1$, which we consider here, $1<\alpha\le 2$.
The endpoint $\alpha=2$ of the interval $\alpha\in(1,2]$ corresponds
to the marginal value $\varepsilon=0$.
Keeping $\varepsilon$ nonnegative means
that we do not move in dimensions $d>4$ where the underlying $\phi^4_d$
theory becomes non-renormalizable.
On the other hand, with the constraint $\alpha> 1$ all Hurwitz and Riemann
zeta functions appearing in the present section are well defined in the
classical sense of (\ref{ZD})--(\ref{e12}).

By construction, the improper integral $K(\alpha)$ is convergent and
well defined (at least) in the whole region of the parameter $\alpha$
indicated in its definition (\ref{e43}), with smooth behavior in the
limit $\alpha\to 2$.
The finite value of $K(\alpha)$ at $\alpha=2$ directly contributes to
the value of the coefficient $c_1(\varepsilon)$ at $\varepsilon=0$ given between
(4.8) and (4.9) in \cite{Sym81}.

The calculation of the integral $K(\alpha)$ defined by (\ref{e43})
is the main nontrivial task of the present section.
In its integrand, we express each Hurwitz zeta function $\zeta(\alpha,y)$
in terms of the auxiliary function $\zeta_1(\alpha,y)$ with
$y\mapsto x$ and $y\mapsto 1-x$. Thus, some straightforward algebra
leads us to
\begin{equation}\label{e46}
K(\alpha)=-4\zeta^2(\alpha)+2I(\alpha,\alpha)+2J(\alpha,\alpha)+2M(\alpha),
\end{equation}
where the integrals $I(\alpha,\alpha)$ and $J(\alpha,\alpha)$
are those from (\ref{e21}) and (\ref{e31}) respectively, and
\begin{align}\label{MA}
&M(\alpha)=
\\\nonumber
&\int_0^1\!\!\Big\{\!\!\!\left[ x^{-\alpha}+(1-x)^{-\alpha}-2\zeta(\alpha) \right]
\left[\zeta_1(\alpha,x)+\zeta_1(\alpha,1-x)-2\zeta(\alpha) \right]
+x^{-\alpha}(1-x)^{-\alpha}\!\Big\}dx.
\end{align}
Note that the integral $M(\alpha)$ is completely symmetric with respect to the transformation
$x\leftrightarrow 1-x$. From the fact that
\begin{align*}
\zeta_1(\alpha,x)&=\zeta(\alpha)-\alpha x\zeta(\alpha+1) +O(x^2),\\
\zeta_1(\alpha,1-x)&=\zeta(\alpha)-1+\alpha x\left[\zeta(\alpha+1)-1\right]+O(x^2),
\end{align*}
it follows that the integrand of $M$ is $O(x^{2-\alpha})+O(1)$ as $x\ra 0$.
By symmetry under the transformation $x\mapsto 1-x$ the behavior of the integrand is
$O((1-x)^{2-\alpha})+O(1)$ as $x\ra 1$.
Hence the integral $M(\alpha)$ converges for $1<\alpha\le 2$.

Substitution of all needed explicit evaluations of Sections 2 and 3
into (\ref{e46}) and (\ref{MA})  then yields
\begin{equation}\label{e47}
K(\alpha)=4\zeta^2(\alpha)+2\zeta(2\alpha-1)\bl\{2B(1-\alpha,2\alpha-1)
+B(1-\alpha,1-\alpha)\br\}+\frac{2}{2\alpha-1}\,.
\end{equation}
Combining (\ref{e40})--(\ref{e42}), (\ref{e44}), and (\ref{e47}) we obtain,
after some straightforward algebra,
the first-order contribution to the Casimir energy (\ref{CE})
\begin{align*}
g (\mu L)^\varepsilon L^{-3+\varepsilon} c_1(\varepsilon)&=
\left\{2^{-9+2\varepsilon}\pi^{-4+\varepsilon}\Gamma^2(1-\varepsilon/2)
g (\mu L)^\varepsilon L^{-3+\varepsilon}\right\}\\
&\times\left\{\zeta^2(\alpha)+(1-\cos\pi\alpha)B(1-\alpha,2\alpha-1)
\zeta(2\alpha-1)\right\}.
\end{align*}
The numerical coefficient of this last equation with $\alpha=2-\varepsilon$
then yields the coefficient $c_1(\varepsilon)$ given by
\begin{equation}\label{e48}
c_1(\varepsilon)=
\frac12 (4\pi)^{-4+\varepsilon}\Gamma^2\Big(1-\frac\varepsilon2\Big)\Big[\zeta^2(2-\varepsilon)+
(1-\cos\pi\varepsilon)B(3-2\varepsilon,-1+\varepsilon)\zeta(3-2\varepsilon)\Big]
\end{equation}
for $\varepsilon\in[0,1)$, as stated in \cite[p.~12]{Sym81}.

As we mentioned before, this result has been obtained by other means
in \cite[Eq. (5.11)]{KD92}, in the context of statistical physics.
There, the constant $c_1(0)=2^{-11}/9$, see \cite[p.~12]{Sym81}, enters the
$O(\varepsilon)$ term of the
$\varepsilon$ expansion of the Casimir amplitude for a Dirichlet film in the
$n$ component $\phi^4_{4-\varepsilon}$ field theory \cite[Eq. (5.16)]{KD92},
\begin{equation}\label{Deld}
\Delta^{\rm DD}=-n\,\frac{\pi^2}{1440}\left[1+\Big(\frac{\gamma-1}{2}
+\ln(2\sqrt\pi)-\frac{\zeta'(4)}{\zeta(4)}-\frac{5}{4}\frac{n+2}{n+8}\Big)
\varepsilon\right]+O(\varepsilon^2)\,.
\end{equation}
Further $\varepsilon$ expansion terms of $c_1(\varepsilon)$ would
contribute to $O(\varepsilon^2)$ and higher order terms that follow
in (\ref{Deld}).
Further generalizations and extensions of this approach
can be found in \cite{Krech,DGS06,BDS10}.

\section{Concluding remarks}

In the present article we have studied the long-standing mathematical problem of calculating
the definite integrals of products of two Hurwitz zeta functions, $\zeta_1(a,x)\zeta_1(b,y)$
with $a,b\in\mathbb C$ and $y=\{x,1-x\}\in[0,1]$.
The problem statement goes back to the classical work of Koksma and Lekkerkerker
\cite{KL52}, where it was formulated as a mean-value theorem for $\zeta_1(s,x)$
with $s=\sigma+it$ and with special emphasis on the behavior of the integral
$\int_0^1|\zeta_1(s,x)|^2 dx$ in the asymptotic limit $t\to\pm\infty$ when $\sigma=\frac12$.
In our notation, this problem statement would correspond to the choice of two
complex conjugate parameters $a=b^*$.
An illuminating survey of this problem and its numerous generalizations can be found in
\cite{Mats00}.

On the other hand, we noticed a close connection of the above mathematical problem
with a possible way \cite{Sym81} of calculating the Feynman integrals
in Euclidean $\phi^4$ field theories
considered in the space between two parallel planes. Here we considered the case of Dirichlet
boundary conditions imposed on boundary planes, though systems with some other boundary
conditions \cite{KD92,Krech} could be treated in a similar manner.
In this context, the case of equal parameters $a=b$ in the zeta functions $\zeta_1$ is of relevance.
We recall that the results derived in Section 4 have been recorded by
Symanzik in \cite{Sym81} without proof and the details of his calculation are not known.
These results were reproduced by Krech and Dietrich
in \cite{KD92} in a different and much less involved way using the momentum representation of the
same Feynman integrals. To the best of our knowledge, our derivation of (\ref{e48}) in terms of integrals
of Hurwitz zeta functions is the first derivation of this formula carried out in this manner.

In doing this, we derived the known result for $I(a,b)$ in a way that
apparently has not been employed by mathematicians for this integral.
Also we needed an analogous result for $J(a,b)$, which has been considered to a lesser
extent in the mathematical literature. Both these integrals were necessary to reconstruct the
coefficient $c_1(\varepsilon)$ from \cite[p. 12]{Sym81} in the main result
(\ref{e48}) of Section 4.
Our motivation for the calculations of Section 4 was to
give a computational scheme for treating more complicated integrals of the same
structure, where the functions in the integrand would be essentially more
complicated but have the same analytical properties.
Such Feynman integrals are needed in an eventual calculation of the next-to-leading,
two-loop contribution to
the Casimir amplitude for strongly anisotropic critical systems in a film geometry
with Dirichlet-Dirichlet boundary conditions. The corresponding one-loop term
has been calculated in \cite[Sec. 5.3.2]{BDS10}.

Finally, we mention that our non-trivial calculations of the integrals considered
in this paper at integer $a=b:=m\ge2\in\mathbb N$ appear to be of both academic and practical interest.

\begin{center}
{\bf Appendix: \ Evaluation of $\bm{I(a,a)}$ for integer $\bm{a\geq 2}$}
\end{center}
\setcounter{section}{1}
\setcounter{equation}{0}
\renewcommand{\theequation}{\Alph{section}.\arabic{equation}}

We write (\ref{e25}) when $a=b$ in the form
\begin{equation}\label{e26}
I(a,a)=\frac{1}{2a-1}-2\,\frac{\Gamma(1-a)}{\Gamma(a)}\Big\{S(a)-\Gamma(2 a-1)\zeta(2 a-1)\Big\},
\end{equation}
where
\[S(a):=\sum_{n=0}^\infty\frac{\g(n+a)}{\g(n+2-a)}\,\{\zeta(a+n)-1\}\]
and the common factor $-2\g(1-a)/\g(a)$ has been extracted from the expression in braces. Use of the well-known result
$\g(n+\alpha)/\g(n+\beta)\sim n^{\alpha-\beta}$ ($n\ra\infty$) shows that the sum $S(a)$ converges absolutely since
the terms are controlled by $n^{2a-2}2^{-n}$ as $n\ra\infty$.

Considering the sum $S(a)$ for integer $a=m\ge2$, we observe that the  function $\g(n+2-m)$ in the denominator
of $S(m)$ has non-positive integer arguments $n+2-m\le0$ and so becomes singular when $n\le m-2$.
This means that the corresponding part of the sum $S(m)$ vanishes, and
the non-vanishing contributions to $S(m)$ come from the terms with
$n\geq m-1$. Thus, upon application of (\ref{e13}) we have
\[
S(m)=\sum_{n=m-1}^\infty \frac{\g(n+m)}{\g(n+2-m)}\,\{\zeta(m+n)-1\}=
\g(2m-1)\sum_{k=0}^\infty\frac{(2m-1)_k}{k!}\,\zeta(2m-1+k,2).
\]
Hence, by (\ref{TE}) with $a\mapsto 2m-1$, $b=2$ and $z=1$
we find
\begin{equation}\label{e27}
S(m)=\g(2m-1)\zeta(2m-1).
\end{equation}
Thus, the difference of terms in braces in (\ref{e26}) vanishes
and the right-hand side of (\ref{e26}) is consequently regular at integer $a=m\ge2$.

We now set $a=m+\epsilon$, where $m=2, 3, \ldots\,$, in $S(a)$ and consider the limit $\epsilon\ra 0$ using the small-parameter expansions (\ref{GE})
and (\ref{EG}). At the moment we have $\alpha=-\epsilon$ in (\ref{EG})
and take only its pole part into account.
By writing $S(a)$ as
\[
S(m+\epsilon)=\sum_{n=0}^{m-2}\frac{\g(n+m+\epsilon)}{\g(-p-\epsilon)}\{\zeta(n+m+\epsilon)-1\}+
\!\sum_{n=m-1}^\infty\!\!\frac{\g(n+m+\epsilon)}{\g(q-\epsilon)}\{\zeta(n+m+\epsilon)-1\}
\]
with $p:=m-n-2\ge0$ and $q:=n+2-m\ge1$ and expanding to first order in $\epsilon$, we find
\begin{align*}
&S(a)=S(m)-\epsilon\sum_{n=0}^{m-2}(-)^{m-n}(n+m-1)!\,(m-n-2)!\,\{\zeta(n+m)-1\}\\
&\hspace{5cm}+\;\epsilon\sum_{n=m-1}^\infty \frac{(n+m-1)!}{(n-m+1)!}\,\Upsilon_n(m)+O(\epsilon^2),
\end{align*}
where we have defined
\begin{equation}\label{e27a}
\Upsilon_n(m):=\{\psi(n+m)+\psi(n+2-m)\} \{\zeta(n+m)-1\}+\zeta'(n+m).
\end{equation}

Similarly, for the last term from (\ref{e26}) we have by (\ref{e27})
\[
\Gamma(2 a-1)\zeta(2 a-1)=S(m)+2\epsilon\Gamma(2m-1)\Big\{\!\psi(2m-1)
\zeta(2m-1)+\zeta'(2m-1)\!\Big\}\!+O(\epsilon^2).
\]
Finally, upon noting that the factor in braces in (\ref{e26}) is
\[\frac{\g(1-a)}{\g(a)}=\frac{(-)^m}{\epsilon\,\g^2(m)}\{1+O(\epsilon)\},\]
we then arrive at
\begin{align}\nonumber
I(m,m)&=\frac{1}{2m-1}+\frac{2(-)^m}{\g^2(m)}
\bl\{ 2\g(2m-1) Z(m)
-\sum_{n=m-1}^\infty \frac{(n+m-1)!}{(n-m+1)!}\,\Upsilon_n(m)\br\}
\\\label{a1}
&+\frac{2}{\g^2(m)}\sum_{n=0}^{m-2}(-)^n\,(n+m-1)!\,(m-n-2)!\,\{\zeta(n+m)-1\},
\end{align}
where $Z(m):=\psi(2m-1) \zeta(2m-1)+\zeta'(2m-1)$.

The infinite sum in (\ref{a1}) can be simplified by removing the term $\zeta'(n+m)$ contained in $\Upsilon_n(m)$.
If we differentiate the well-known result obtained from (\ref{e13}) with $b=2$, $z=1$
\[\sum_{k=1}^\infty \frac{\g(\alpha+k)}{k!} \{\zeta(\alpha+k)-1\}=\g(\alpha)\]
with respect to $\alpha$, we find that
\[\sum_{k=1}^\infty \frac{\g(\alpha+k)}{k!}\,\zeta'(\alpha+k)=\g(\alpha)\psi(\alpha)\zeta(\alpha)-\sum_{k=0}^\infty \frac{\g(\alpha+k)}{k!} \,\psi(\alpha+k)\{\zeta(\alpha+k)-1\}.\]
Then, the contribution from the part of the sum containing $\zeta'(n+m)$ in (\ref{a1}) is
\begin{align*}
\sum_{n=m-1}^\infty &\frac{(n+m-1)!}{(n-m+1)!}\,\zeta'(n+m)\\
&=\g(2m-1)
\zeta'(2m-1)+\sum_{k=1}^\infty\frac{\g(2m-1+k)}{k!}\,\zeta'(2m-1+k)\\
&=\g(2m-1) Z(m)-\sum_{k=0}^\infty \frac{\g(2m-1+k)}{k!}\,\psi(2m-1+k)\{\zeta(2m-1+k)-1\}\\
&=\g(2m-1) Z(m)-\sum_{n=m-1}^\infty \frac{(n+m-1)!}{(n-m+1)!}\,\psi(n+m) \{\zeta(n+m)-1\}.
\end{align*}
Upon substitution of this evaluation in (\ref{a1}) we arrive at
the result stated in (\ref{e28}) in Theorem 2.

\providecommand{\bysame}{\leavevmode\hbox to3em{\hrulefill}\thinspace}
\providecommand{\MR}{\relax\ifhmode\unskip\space\fi MR }
\providecommand{\MRhref}[2]{%
  \href{http://www.ams.org/mathscinet-getitem?mr=#1}{#2}
}
\providecommand{\href}[2]{#2}

\end{document}